\documentclass[12pt,a4paper]{amsart} 
\usepackage{mathrsfs}
\usepackage{mathtools, booktabs}
\usepackage{amstext} 
\usepackage{array}   
\newcolumntype{R}{>{$}r<{$}} 
\usepackage{dsfont}

\theoremstyle{plain}                 
\newtheorem{theorem}{Theorem}[section]     
\newtheorem{proposition}[theorem]{Proposition} 
\newtheorem{corollary}[theorem]{Corollary}     
        
\theoremstyle{definition}           
\newtheorem{definition}{Definition}    
 
\theoremstyle{remark}       
\newtheorem{remark}{Remark}

\DeclareMathOperator\Id{Id}
\DeclareMathOperator\Class{Class}
\DeclareMathOperator\NonPoly{NonPoly}
\DeclareMathOperator\End{End}

\newcommand{\sta}{\mathcal}

\newcommand{\MMMbar}{\overline{\sta M}}

\def\cM{{\mathcal{M}}}
\def\cL{{\mathcal{L}}}
\def\cC{{\mathcal{C}}}
\def\Klog{\omega_{\rm log}}
\def\oM{\overline{\mathcal{M}}}
\def\Z{\mathbb{Z}}

\def\ch{\mathrm{ch}}
\def\C{\mathrm{C}}
\def\H{\mathrm{H}}
\def\L{\mathrm{L}}
\def\S{\mathcal{S}}
\def\bbC{\mathbb{C}}
\def\oM{\overline{\mathcal{M}}}
\def\E{\mathrm{E}}
\def\V{\mathrm{V}}

\def\Aut{{\rm Aut}}

\def\dd{\mathrm{d}}
\def\CP1{\mathbb{C}\mathrm{P}^1}

\makeatletter
\renewenvironment{cases}[1][l]{\matrix@check\cases\env@cases{#1}}{\endarray\right\rbrace}
\def\env@cases#1{%
  \let\@ifnextchar\new@ifnextchar
  \left.\def\arraystretch{1.2}%
  \array{@{}#1@{\, }l@{}}
  }
\makeatother

\begin{document}
\title[On ELSV, Hurwitz numbers and topological recursion]{On ELSV-type formulae, Hurwitz numbers and topological recursion}

\author[D.~Lewanski]{D.~Lewanski}
\address{D.~L.: Korteweg-de Vries Institute for Mathematics, University of Amsterdam, Postbus 94248, 1090 GE Amsterdam, The Netherlands}
\email{D.Lewanski@uva.nl}

\begin{abstract}
We present several recent developments on ELSV-type formulae and topological recursion concerning Chiodo classes and several kind of Hurwitz numbers. The main results appeared in~\cite{LPSZ}. 
\end{abstract}

\maketitle

\tableofcontents

\section{Introduction}
ELSV-type formulae relate connected Hurwitz numbers to the intersection theory of certain classes on the moduli spaces of curves. Both Hurwitz theory and the theory of moduli spaces of curves benefit from them, since ELSV formulae provide a bridge through which calculations and results can be transferred from one to the other. The original ELSV formula \cite{ELSV} relates simple connected Hurwitz numbers and Hodge integrals. 
It plays a central role in many of the alternative proofs of Witten's conjecture that appeared after the first proof by Kontsevich (for more details see \cite{LIU}).

\subsubsection{Examples of ELSV-type formulae:}
The simple connected Hurwitz numbers $h^{\circ}_{g;\vec{\mu}}$ enumerate connected Hurwitz coverings of the $2$-sphere of degree $|\vec{\mu}|$ and genus $g$, where the partition $\vec{\mu}$ determines the ramification profile over zero, and all other ramifications are simple, i.e. specified by a transposition of two sheets of the covering $(a_i \, b_i) \in \mathfrak{S}_{|\vec{\mu}|}$. By the Riemann Hurwitz formula, the number of these simple ramifications is $b = 2g - 2 + l(\vec{\mu}) + |\vec{\mu}|$ (for an introduction on Hurwitz theory see, e.g., \cite{CM}).\\
The celebrated ELSV formula \cite{ELSV} expresses these numbers in terms of the intersection theory of moduli spaces of curves:
\begin{align*}
 \frac{ h^{\circ}_{g;\vec{\mu}} }{b!} = 
  \prod_{i=1}^{\ell(\vec{\mu})}\frac{\mu_i^{\mu_i}}{\mu_i !}
  \int_{\oM_{g,\ell(\vec{\mu})}}
  \left(\sum_{l=0}^g(-1)^l \lambda_l \right) \prod_{j=1}^{\ell(\vec{\mu})} \sum_{d_j = 0}\mu_j^{d_j} {\psi}_j^{d_j}
 \end{align*}
A different Hurwitz problem arises from Harish-Chandra-Itzykson-Zuber matrix model and leads to the study of the simple monotone connected Hurwitz numbers $h_{g,\mu}^{\circ, \le}$ (see \cite{GGN}), in which an extra \textit{monotone} condition is imposed on the coverings \textemdash \, if $(a_i \, b_i)_{i=1, \dots, b}$ are written such that $a_i < b_i$, the condition requires that $b_i \geq b_{i+1}$ for all $i = 1, \dots, b-1$. For these Hurwitz numbers an ELSV-type formula is also known \cite{ALS, DK}:
	\begin{align*} 
	h_{g,\vec{\mu}}^{\circ, \le}&  =  \prod_{i=1}^{\ell(\vec{\mu})} \binom{2\mu_i}{\mu_i} 
	\int_{\overline{\mathcal{M}}_{g,\ell(\vec{\mu})}} 
	 \!\!\!\!\!\! \exp\left(\sum_{l=1} A_l \kappa_l \right)
	\prod_{j=1}^{\ell(\vec{\mu})} \sum_{d_j = 0} \psi_j^{d_j} \frac{(2(\mu_j + d_j) - 1)!!}{(2 \mu_j - 1)!!},
	\end{align*}
where the generating series for the coefficients $A_i$ of the kappa classes reads
$
	\exp \left(-\sum_{l=1}^\infty A_l U^l \right)  = \sum_{k=0}^{\infty} (2k + 1)!! U^k.
$
More examples exist in the literature, and two of them will be subject of this paper: the ELSV formula due to Johnson, Pandharipande and Tseng (JPT) \cite{JPT} and a conjectural formula due to Zvonkine \cite{Z}, see also \cite{SSZ} (for the precise formulae see the table at the end of this paper).\\
\subsubsection{Structure of ELSV-type formulae:}
All four examples above express numbers enumerating connected Hurwitz covers of a certain kind, depending on a genus parameter and a partition, in terms of some \textit{non - polynomial} factor in the entries of a partition $\vec{\mu}$ (in the formulae above $\prod \frac{{\mu_i}^{\mu_i}}{\mu_i !}$ and $ \binom{2\mu_i}{\mu_i} $) and an integral over moduli spaces of curves of a certain class intersected with $\psi$ class. This integral is clearly a \textit{polynomial} of degree $3g - 3 + \ell(\vec{\mu})$ in the $\mu_i$. Conceptually:
\begin{equation}\label{eq:structureELSV}
h^{\circ, condition}_{g,\vec{\mu}} = \NonPoly (\vec{\mu}) 
\int_{\overline{\mathcal{M}}_{g,\ell(\vec{\mu})}}
\!\!\!\!\!\!\!\! (\Class) \,\prod_{j=1}^{\ell(\vec{\mu})} \sum_{d_j=0} c_{d_j}(\mu_j) \psi_j^{d_j}
\end{equation}
where $c_{d_j}(\mu_j)$ is a polynomial of degree $d_j$ in $\mu_j$.
\subsubsection{Topological recursion and ELSV formulae}
The Chekhov, Eynard and Orantin (CEO) topological recursion procedure associates to a spectral curve $\S = (\Sigma, x(z), y(z), B(z_1, z_2))$ (see, e.g., \cite{CE, EO}) a collection of symmetric correlation differentials $\omega_{g,n}$ defined on the product of the curve $\Sigma^{\times n}$ through a universal recursion on $2g - 2 + n$. The expansion of these differentials near particular points can unveil interesting invariants, or solutions to enumerative geometric problems.\\
\begin{center}
\includegraphics[width=11cm]{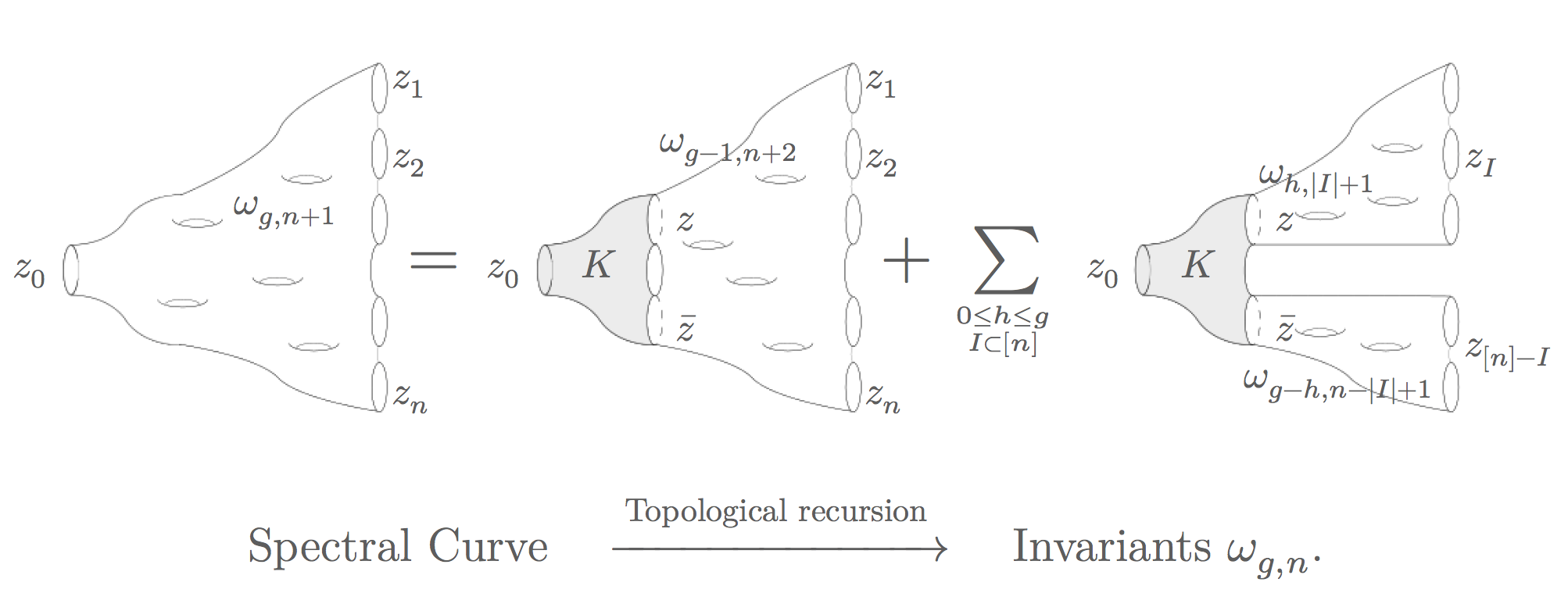}
\end{center}
We say that certain numbers satisfy the topological recursion if there exists a spectral curve such that the expansion of the correlation differentials near \textit{some} point has those numbers as coefficients.
The expansion of the correlation differentials takes the form
\begin{equation}
\omega_{g,n}^{\S} = \dd_1 \otimes \dots \otimes \dd_n \sum_{\mu_1, \dots, \mu_n} N^{\S}_{g,\vec{\mu}} \prod_{i=1}^n\tilde{x}_i^{\mu_i}
\end{equation}
for some coefficients $N_{g,\vec{\mu}}$, where $\tilde{x}$ is a function of $x$ that depend on the point of the expansion.
Both the simple Hurwitz and the monotone Hurwitz numbers satisfy the topological recursion (see \cite{BEMS, BM, DDM, EMS, MZ}), and their spectral curves are respectively
\begin{equation}\label{eq:2spetral}
\left(\mathbb{C} \mathbb{P}^1, -z + \log(z), z, \frac{\dd z_1 \dd z_2}{(z_1 - z_2)^2} \right), \left(\mathbb{C} \mathbb{P}^1, \frac{z-1}{z^2}, -z, \frac{\dd z_1 \dd z_2}{(z_1 - z_2)^2}\right)
\end{equation}
In the simple Hurwitz case $\tilde{x} = e^x$, whereas in the monotone case $\tilde{x} = x$. On the other side, it was proved that the expansion of the correlation differentials \textit{have the same structure as the right-hand side of ELSV-type formulae} described above (see Theorem \ref{thm:CohFT}), depending on same ingredients that are functions of the spectral curve.\\
At this point comes the key observation: if one can compute these ingredients explicitly for a given spectral curve $\S$, one proves that 
\begin{equation}\label{eq:structureELSV}
N_{g, \vec{\mu}}^{\S} = \NonPoly^{S} (\vec{\mu}) 
\int_{\overline{\mathcal{M}}_{g,\ell(\vec{\mu})}}
\!\!\!\!\!\!\!\! (\Class^{\S}) \,\prod_{j=1}^{\ell(\vec{\mu})} \sum_{d_j=0} c_{d_j}^{\S}(\mu_j) \psi_j^{d_j},
\end{equation}
where the non-polynomial part $\NonPoly^{\S}$, the class $Class^{\S}$, and the $c_{d_j}^{\S}(\mu_j)$ are explicit. This allows to formulate equivalence statements in the following sense. 
\begin{definition}\label{def:eqst}
An \textit{TR-ELSV equivalence statement} for a Hurwitz problem $h_{g, \vec{\mu}}^{\circ, condition}$ and a spectral curve $\S$ asserts the equivalence between the following two propositions:
\begin{enumerate}
\item[i)] The numbers $h_{g, \vec{\mu}}^{\circ, condition}$ satisfy the topological recursion with input spectral\ curve $\S$ (i.e. $h_{g, \vec{\mu}}^{\circ, condition} = N_{g, \vec{\mu}}^{\S}$) \\
\item[]
$$\text{ii)}\,\,
h^{\circ, condition}_{g,\vec{\mu}} = \NonPoly^{\S} (\vec{\mu}) 
\int_{\overline{\mathcal{M}}_{g,\ell(\vec{\mu})}}
\!\!\!\!\!\!\!\! (\Class^{\S}) \,\prod_{j=1}^{\ell(\vec{\mu})} \sum_{d_j=0} c_{d_j}^{\S}(\mu_j) \psi_j^{d_j}.$$
\end{enumerate}
\end{definition}

Clearly, it makes sense to formulate an equivalence statement for certain Hurwitz numbers $h^{\circ, condition}_{g,\vec{\mu}}$ and a certain spectral curve $\S$ if for at least one of the two propositions there exists some evidence or a proof. \\

Thus, if one establishes $i)$ independently of $ii)$, then $ii)$ follows immediately (and vice versa), and hence this equivalence relationship has received much attention in the literature. 
For example, for the case of simple Hurwitz numbers and the first curve in Equation \eqref{eq:2spetral}, proposition $i)$ was conjectured to hold by Bouchard and Mari\~{n}o \cite{BM}, while $ii)$ is the original ELSV formula. Proposition $i)$ was proved in \cite{BEMS, EMS, MZ}, the equivalence statement was proved in \cite{Eyn11}, see also \cite{SSZ}. The equivalence statement immediately provides a new proof of $i)$ from $ii)$. The proofs \cite{EMS, MZ} of $i)$ though, make use of a polynomiality property that is extracted from ELSV, hence the equivalence cannot be used in the other direction without falling into a circular argument, unless this polynomiality property can be proved without using ELSV formula. 
This was done in \cite{DKOSS}, see also \cite{DLPS, KLS}, 
 and thus $ii)$ follows from $i)$ by the equivalence statement.\\
In the case of r-spin Hurwitz numbers, proposition $i)$ is known as $r$-Bouchard-Mari\~{n}o conjecture \cite{BM}, proposition $ii)$ is the $r$-ELSV formula conjectured by Zvonkine \cite{Z}, see also \cite{SSZ}. The equivalence of $i)$ and $ii)$ was established in \cite{SSZ}, but since neither of the two have been proven, both remain conjectural.


%
\subsubsection{Givental theory}
The class $\Class^{\S}$ describes a semi-simple cohomological field theory (CohFT), possibly with non-flat unit. Semi-simple CohFTs with unit are classified by Givental - Teleman \cite{G,T} by the action of a Givental $R$-matrix on a topological field theory. On the other hand, Givental theory has been identified with the CEO topological recursion in \cite{DOSS}. This identification makes Givental theory a powerful and explicit tool to compute the ingredients above, and hence to prove equivalence statements.

\subsubsection{Main result}
The main result in \cite{LPSZ} is the computation of the ingredients $\NonPoly^{\S_{r,s}}, \, \Class^{\S_{r,s}}, $ and $c_{d_j}^{\S_{r,s}}(\mu_j)$ for the specific spectral curve 
\begin{equation*}
	\S_{r,s} := \left( \CP1, x(z)=-z^r+\log z,
 y(z)=z^s,  B(z,z') = \frac{\dd z\, \dd z'}{(z-z')^2} \right).
\end{equation*} 
The class $\Class^{\S_{r,s}}$ turns out to coincide with Chiodo classes \cite{Chiodo}. 
Its specialisation for $r=s=1$ is used in the ELSV-type formulae for simple Hurwitz numbers, the case $s=1$ involves the $r$-spin Hurwitz numbers, while the case $r=s$ is used in the ELSV formula for $r$-orbifold Hurwitz numbers, derived by Johnson, Pandharipande and Tseng (JPT) \cite{JPT}. The equivalence statement for general $r$ and $s$ is derived, and it specialises to the equivalence statements already known for simple and $r$-spin Hurwitz numbers. For $r$-orbifold Hurwitz instead, the statement is new (see Corollary \ref{cor:eqorbifold2}): this gives a new proof of the topological recursion for $r$-orbifold Hurwitz numbers from JPT. On the other hand, the topological recursion was already proved in \cite{BHLM, DLN}, but extracting some polynomiality property from JPT formula itself. This polynomiality property was then proved in \cite{DLPS}, see also \cite{KLS}, providing together with the equivalence statement a new proof of JPT formula.\\

\subsubsection{Plan and organisation of the paper}
In this note we give a short exposition of the results obtained in \cite{LPSZ}. On the one hand, the technical proofs are omitted and we refer to the original paper for them; conversely, some of the details that are omitted in the paper are here worked out. \\
In Section \ref{sec:Chiodo} we recall the definition of Chiodo classes, we express them in terms of stable graphs and prove that they are given by the action of a particular Givental $R$-matrix.\\
In Section \ref{sec:CurvetoGivental} we review the key steps of DOSS Identification and we show the result of the computation for the spectral curve $\S_{r,s}$. This allows to state the main result (see Theorem \ref{TH2}).\\
In Section \ref{sec:JPT} we treat the equivalence statement for $r$-orbifold Hurwitz numbers and Johnson-Pandharipande-Tseng formula.

\subsection{Acknowledgments} 
This note is based on a talk given by the author at the 
2016 AMS von Neumann Symposium
\textit{Topological Recursion and its Influence in Analysis, Geometry, and Topology}, July 4 \textemdash \, 8, 2016, in Charlotte, NC.
I would like to thank the organisers of the Symposium and the AMS for this opportunity.\\
I would moreover like to thank N.~Do, P.~Dunin-Barkowski, M.~Karev, R.~Kramer, A.~Popolitov, S.~Shadrin and D.~Zvonkine for interesting discussions, special thanks to S.~Shadrin for having introduced me to the topic and very useful remarks. The author is supported by the Netherlands Organization for Scientific Research.

\section{Chiodo classes}\label{sec:Chiodo}
In this section we recall Chiodo classes and we show their Givental decomposition, for more details we refer the reader to \cite{Chiodo, ChiRua, ChiZvo, JPPZ, SSZ}.
For $2g - 2 +n > 0$, consider a nonsingular curve with distinct markings $[C, p_1, \dots, p_n]\in \cM_{g,n}$, and let $\Klog = \omega_C(\sum p_i )$ be its log canonical bundle.
Let $r \geq 1$, $1\leq a_1,\dots, a_n \leq r$ and $0 \leq s \leq r$ be integers satisfying the condition
\begin{equation}\label{eq:cond_a}
(2g-2+n)s-\sum_{i=1}^n a_i = 0 \quad \mod \, r
\end{equation}
This condition guarantees the existence of $r$th tensor roots $L$ of the line bundle
 $\Klog^{\otimes s}\left(-\sum 
a_i p_i\right)$ on $C$. 
For the moduli space of such $r$th tensor roots, a natural compactification $\oM_{g;a_1, \dots, a_n}^{r,s}$ was constructed in~\cite{Chiodocostruzione, Jarvis}. Let
$\pi : \cC_{g;a_1, \dots, a_n}^{r,s} \to \oM_{g;a_1, \dots, a_n}^{r,s}$
 be the universal curve, let $\cL \to \cC_{g;a_1, \dots, a_n}^{r,s}$ be the universal 
$r$th root, and let
$ \epsilon: \oM_{g;a_1, \dots, a_n}^{r,s} \to \oM_{g,n}$
 be the forgetful map (in order for $\epsilon$ to be unramified in the orbifold sense, the target $\oM_{g,n}$ is changed into the moduli space $r$-stable curves, meaning that for each stable curve there is an extra $\Z_r$ stabilizer at each node, see \cite{Chiodocostruzione}).
Recall the generating series for the Bernoulli polynomials
$$\sum_{l=0} B_l(x) \frac{t^l}{l!} = \frac{t e^{xt}}{e^t  - 1}$$
where the usual Bernoulli numbers are $B_l(0) = B_l$.
We are interested in the Chiodo classes \cite{Chiodo}
\begin{align*}\label{eq:DefinitionChiodoCohFT}
& \C_{g,n}(r,s;a_1,\dots,a_n):=   \epsilon_{*} c\big(-R^*\pi_*\cL\big) =\\ \notag
& \epsilon_{*}\exp\left(\sum_{l=1}^\infty (-1)^l (l-1)!\ch_l(r,s;a_1,\dots,a_n)\right) \in H^{even}(\oM_{g,n}),
\end{align*}
where Chiodo formula for the Chern characters reads
\begin{align}
\ch_l(r,s;a_1,\dots,a_n)  =		
		\frac{B_{l+1}(\frac sr)}{(l+1)!} \kappa_l
		- \sum_{i=1}^n 
		\frac{B_{l+1}(\frac{a_i}r)}{(l+1)!} \psi_i^l 
		\\ \notag 
		+ \frac{r}2 \sum_{a=1}^{r} 
		\frac{B_{l+1}(\frac{a}r)}{(l+1)!} (j_a)_* 
		\frac{(\psi')^l + (-1)^{l-1} (\psi'')^l}{\psi'+\psi''}.
\end{align}
Here $j_a$ is the boundary map that represents the boundary divisor with remainder $a$ at one of the two half edges, and $\psi',\psi''$ are the $\psi$-classes at the two branches of the node.\\
For the specialisation $r = s =1$, and moreover $a_i = 1$, for $i = 1, \dots, n$, the map $\epsilon$ is the identity map and we recover Mumford formula \cite{Mu_GRR} for the total Chern class of the dual of the Hodge bundle $c(\Lambda^{\vee}_g)$:
\begin{align}
C_{g,n}(1,1; 1, \dots, 1) &=\exp \Bigg( - \Bigg[ \sum_{l=1}^{\infty} 
\frac{B_{l+1}}{l(l+1)} \kappa_l
		- \sum_{i=1}^n 
		\frac{B_{l+1}}{l(l+1)} \psi_i^l 
		\\ \notag
		&+ \frac{1}2 
	    \frac{B_{l+1}}{l(l+1)} j_* 
		\frac{(\psi')^l + (-1)^{l-1} (\psi'')^l}{\psi'+\psi''} \Bigg]\Bigg)\\ \notag
		&= c(\Lambda^{\vee}_g ) = 1 - \lambda_1 + \lambda_2 - \dots + (-1)^g \lambda_g 
\end{align}
where the identity $B_l(1) = (-1)^l B_l$ is used. The formula in \cite{Mu_GRR} is slightly different due to a different Bernoulli number convention and a missprint in the $\kappa$ term.

\subsection{Expression in terms of stable graphs}
Let us recall the expression of the Chiodo class in terms of the sum of products of contributions decorating stable graphs, in order to compare it with the Givental action, for more details see \cite{JPPZ}. The strata of the moduli space of curves correspond
to stable graphs 
$$\Gamma=(\V,\E,\H,\L, g,n :\V \rightarrow \Z_{\geq 0}, v:\H \rightarrow \V, \iota : \H \rightarrow \H)$$
where $\V (\Gamma),\, \E(\Gamma),\, \H(\Gamma)$ and $\L(\Gamma)$ respectively denote the sets of vertices, edges, half-edges and leaves of $\Gamma$; self-edges are permitted. A half-edge indicates either a leaf or an edge together with a choice of one of the two vertices it is attached to. The function $v$ associates to each half-edge its vertex assignment, while $\iota$ is the involution that swaps the two half-edges of the same edge, or leaves the half-edge invariant if it is a leaf. The function $n(v)$ denotes the valence of $\Gamma$ at $v$, including 
both half-edges and legs, and $g(v)$ denotes the genus function. Every vertex $v$ is required to satisfy the stability condition $2g(v) - 2 + n(v) >0$, and the genus of a stable graph $\Gamma$ is defined by $g(\Gamma):= \sum_{v\in V} g(v) + h^1(\Gamma)$.
Let $\Aut(\Gamma)$ denote the group of automorphisms
of the sets $\V$ and $\H$ which leave the
structures $\L$, $\mathrm{g}$, $v$, and $\iota$ invariant.
Let $\mathsf{G}_{g,n}$ be the finite set of isomorphism classes of stable graphs of genus $g$ with $n$ legs. 
Let moreover $\mathsf{W}_{\Gamma,r,s,\vec{a}}$ be the set of \textit{weightings} $ w:\H(\Gamma) \rightarrow \{ 0,\ldots, r-1\}$
satisfying the following three properties:
\begin{enumerate}
\item[(i)] The $i$-th leaf $l_i$ has weight $w(l_i)=a_i \mod r \,$, for $ i\in \{1,\ldots, n\}$.
\item[(ii)] For any two half-edges $h'$ and $h''$ corresponding to the same edge, we have $w(h')+w(h'')=0 \mod r\,$.
\item[(iii)] The condition in Equation \eqref{eq:cond_a} is satysfied locally on each component: for any vertex $v$ the sum of the weights associated to the half-edges incident to $v$ is
$\sum_{v(h)= v} w(h) = s\big( 2g(v)-2+n(v)\big) \mod r $. 
\end{enumerate}

\begin{proposition}\cite{JPPZ}. \label{Cor:ChiodoExp}
The Chiodo class $\C_{g,n}(r,s;a_1,\dots,a_n)\in R^*(\oM_{g,n})$ is equal to 
\begin{multline}\label{eq:ChiodoExp}
\hspace{-10pt}\sum_{\Gamma\in \mathsf{G}_{g,n}} 
\sum_{w\in \mathsf{W}_{\Gamma,r,s,\vec{a}}}
\frac{r^{|E(\Gamma)| + \sum_{v \in V(\Gamma)} 2g(v)-1}}{|\Aut(\Gamma)| }
\;
\xi_{\Gamma*}\Bigg[ \prod_{v \in \V(\Gamma)} e^{-\sum\limits_{l\geq 1} (-1)^{l-1}\frac{B_{l+1}(s/r)}{l(l+1)}\kappa_l(v)} \; \cdot 
\\ 
\prod_{i=1}^n e^{\sum\limits_{l\geq 1}(-1)^{l-1} \frac{B_{l+1}(a_i/r)}{l(l+1)} \psi^l_{h_i}} \cdot 
\prod_{\substack{e\in \E(\Gamma) \\ e = (h',h'')}}
\frac{1-e^{\sum\limits_{l \geq 1} (-1)^{l-1} \frac{B_{l+1}(w(h)/r)}{l(l+1)} [(\psi_{h'})^l-(-\psi_{h''})^l]}}{\psi_{h'} + \psi_{h''}} \Bigg]\, .
\end{multline} 
where $\xi_{\Gamma}$ is the
canonical
morphism 
$
\xi_{\Gamma}: \prod_{v\in \V(\Gamma)} \oM_{g(v),n(v)} \rightarrow \oM_{g,n}
$
 of the boundary stratum corresponding to $\Gamma$.
\end{proposition}

\subsection{Expression in terms of Givental action}\label{sec:giv}
In this section we express Chiodo classes in terms of Givental theory. \\

Fix a vector space $V$ and a symmetric bilinear form $\eta$ on $V$. A Givental $R$-matrix is a $\End(V)$-valued power series
\begin{equation}
R(\zeta) = 1 + \sum_{l=1} R_l \zeta^l = \exp\left(\sum_{l=1}r_l \zeta^l \right), \quad R_l, r_l \in \End(V)
\end{equation}
satisfying the symplectic condition 
$$R(\zeta)R^{*}(-\zeta) = 1 \in \End(V) $$
where $R^*$ is the adjoint of $R$ with respect to $\eta$. 

By Givental - Teleman classification \cite{G,T}, every semi-simple cohomological field theory (CohFT) with unit is obtained by the action of a Givental $R$-matrix on a topological field theory. 
We will show that the action of the $R$-matrix
\begin{align}
	& R^{-1}(\zeta) := \exp\left(-\sum_{l=1}^{\infty} \frac{\mathrm{diag}_{a=1}^{r} B_{l+1}\left(\frac{a}{r}\right)}{l(l+1)}(-\zeta)^l\right)
\end{align}
defined as power series valued in the endomorphisms for the vector space $$V = \langle v_1, \dots, v_{r} \rangle$$ with 
$$\eta(v_a,v_b)= \frac{1}{r}\delta_{a+b \mod r},$$
acting on the topological field theory 
\begin{equation}	
\alpha^{top}_{g,n}(v_{a_1}\otimes \cdots \otimes v_{a_n}) =r^{2g-1} \delta_{a_1+\cdots+a_n-s(2g-2+n) \mod r},
\end{equation}
produces the Chiodo classes. Therefore Chiodo classes determine a CohFT with a known Givental decomposition.
The action of the Givental $R$-matrix is defined as sum over stable graphs $\Gamma$ weighted by ${|\Aut (\Gamma)|}^{-1}$, with contributions on the leaves, on the edges, on special leaves called \textit{dilaton} leaves, and the topological field theory contributes on the verteces. Chiodo classes are already expressed as a sum over stable graphs in Equation \eqref{eq:ChiodoExp} with a very similar structure.
Let us match the Givental contributions one by one:

\subsubsection{Ordinary leaf contributions.} The contribution of the $i$-th leaf reads
$$\exp\Bigg(-\sum_{l=1} \frac{B_{l+1}(a_i/r)}{l(l+1)} (-\psi_{h_i})^{l}\Bigg) = \sum_{j=1}^r (R^{-1})_{a_i}^j(\psi_{h_i})$$
\subsubsection{Dilaton leaf contributions.} Recall that the kappa classes are defined as $\kappa_l = \pi_*(\psi_{n+1}^{l+1})$ under the map that forgets the last marked point $\pi: \oM_{g, n+1} \rightarrow \oM_{g,n}$. 
The contributions on the dilaton leaves correspond to the contributions on the vertices in Equation \eqref{eq:ChiodoExp} before forgetting the corresponding marked point. For the dilaton leaf marked with label $n+i$, for some positive integer $i$, the contribution reads:
$$\exp\Bigg(-\sum_{l=1} \frac{B_{l+1}(s/r)}{l(l+1)} (-\psi_{n+i})^{l} (- \psi_{n+i})\Bigg) $$
We check that $v_s$ is the neutral element $\mathds{1}$ for the quantum product $\bullet$ in flat basis:
\begin{align}
\eta(v_s \bullet v_a, v_b )  &= \alpha^{top}_{0,3}(v_s \otimes v_a \otimes v_b)\\ \notag
 &= r^{-1} \delta_{s + a + b -s  \mod r} = r^{-1} \delta_{ a + b  \mod r} \\ \notag
 & = \eta(v_a, v_b)
\end{align}
Hence the contribution of the dilaton leaf $n+i$ is
$$\psi_{n+i}\Big[ \Id - \sum_{j=1}^r (R^{-1})_{\mathds{1}}^j(\psi_{n+i})\Big]$$
\subsubsection{Edge contributions.} The edge contribution in Equation \eqref{eq:ChiodoExp}, multiplied by the factor
 $(\psi_{h'} + \psi_{h''})$ and after applying the property of Bernoulli numbers
 $(-1)^{p+1} B_{p+1}\left(\frac{w(h')}{r}\right) = B_{p+1}\left(\frac{r-w(h')}{r}\right)$, reads
 $$1 -\, \exp\Bigg( - \sum_{l=1} \frac{B_{l+1}(\frac{w(h')}{r})}{l(l+1)}(-\psi_{h'})^l\Bigg)\exp\Bigg( - \sum_{p=1} \frac{B_{p+1}(\frac{r- w(h'')}{r})}{p(p+1)}(-\psi_{h'})^p\Bigg).$$
Note that the condition on the weightings $w(h') + w(h'') = 0 \mod r$ can be taken care by the scalar product $\eta$. Hence we can write the Givental contribution on the edges as
\begin{equation*}
\sum_{j_1,j_2} \frac{\eta^{j_1,j_2} -  (R^{-1})_{w(h')}^{j_1} (\psi_{h'})\eta^{w(h'), w(h'')}(R^{-1})_{w(h'')}^{j_2}(\psi_{h''})}{\psi_{h'} + \psi_{h''}}
\end{equation*}
\subsubsection{Weightings.} Out of the three conditions on the weightings, condition $(i)$ becomes $w(l_i) = a_i$, condition $(ii)$ on the edges is taken care by the bilinear form $\eta$, condition $(iii)$ can be substituted by the topological field theory condition.
\subsubsection{Powers of $r$.} Every stable graph contributes with
$$
|E(\Gamma)| + \sum_{v \in V(\Gamma)} 2g(v) - 1
$$
powers of $r$. Indeed the topological field theory in the vertex $v$ provides $2g(v) - 1$ powers of $r$, and the inverse of $\eta$ provides one power of $r$ for each edge of $\Gamma$.
\subsubsection{The expression of Givental action}
Let us indicate with $$ \{l_1, \dots, l_n, l_{n+1}, \dots, l_{n+k} \} = L(\Gamma) $$ the set of legs, corresponding to marked points of the curves in $\oM_{g, n+k}$, and let
$$\xi^{(k)}_{\Gamma}: \prod_{v \in V(\Gamma)} \oM_{g(v), n(v)} \rightarrow \oM_{g,n}$$ be the canonical morphism of the boundary stratum corresponding to $\Gamma$ that forgets the last $k$ marked points. Let us consider functions $w^{\vee}: H(\Gamma) \rightarrow \Z_{\geq 0} $ without \textit{any} further condition. We use here the notation $w^{\vee}$, instead of $w$, to remark that the weightings $w^{\vee}$ decorates the half-edges \textit{after} the application of the endomorphisms $R^{-1}_l$. Collecting the contributions and the considerations above, we have:

\begin{align*}
C_{g,n}&(r,s; a_1, \dots, a_n) =\sum_{k = 0} \, \sum_{\substack{\Gamma \in G_{g,n+k}\\ w^{\vee}: H(\Gamma) \rightarrow \Z_{\geq 0}}}  \frac{1}{|\Aut(\Gamma)|} \left(\xi^{(k)}_{\Gamma}\right)_*\Bigg[\\ \notag
&\prod_{v \in V(\Gamma)} \alpha^{top}_{g(v),n(v)}\Bigg(\bigotimes_{\substack{h \in H(\Gamma) :\\ v(h) =v }} v_{w^{\vee}(h)} \Bigg)  \\ \notag
&\prod_{i=1}^n (R^{-1})_{a_i}^{w^{\vee}(l_i)}(\psi_{i})
\\ \notag &
\prod_{i=1}^{k} \psi_{n+i}\Big[ \Id -  (R^{-1})_{\mathds{1}}^{w^{\vee}(l_{n+i})}(\psi_{n+i})\Big]
\\ \notag
&\prod_{e=(h',h'') \in V(\Gamma)}
\!\!\!\!\!\!\!\!\!\!\!\!\ \frac{\eta^{w^{\vee}(h'),w^{\vee}(h'') } - \sum_{k_1, k_2}(R^{-1})^{w^{\vee}(h')}_{k_1} (\psi_{h'})\eta^{k_1, k_2}(R^{-1})^{w^{\vee}(h'')}_{k_2}(\psi_{h''})}{\psi_{h'} + \psi_{h''}}
\Bigg]  \\ \notag
\end{align*}
 The expression above is equivalent to $\left(R. \alpha^{top}\right)_{g,n}(v_{a_1} \otimes \dots \otimes v_{a_n})$, i.e. the Givental action of the matrix $R$ on $\alpha^{top}$ correspondent of genus $g$ and $n$ marked points, evaluated in the element $v_{a_1} \otimes \dots \otimes v_{a_n}$, (see \cite{G, DOSS, PPZ}).\\
Consider then $\C_{g,n}(r,s;a_1,\dots,a_n)$ as the evaluation of a map 
\begin{equation}
\C_{g,n}(r,s)\colon V^{\otimes n}\to H^{even}(\oM_{g,n}),
\end{equation}
where $V=\langle v_1,\dots,v_r\rangle$, and 
\begin{equation}
\C_{g,n}(r,s)\colon v_{a_1}\otimes\cdots\otimes v_{a_n} \mapsto
\C_{g,n}(r,s;a_1,\dots,a_n).
\end{equation}
The previous calculation shows
\begin{proposition}[\cite{LPSZ}]\label{prop:RChiodo}
 For $0\leq s\leq r$ the collection of maps $\{\C_{g,n}(r,s)\}$ defined by the Chiodo classes form a semi-simple cohomological field theory with flat unit, obtained by the action of the Givental matrix $R$ on the topological field theory
 $\alpha^{top}_{g,n}$:
 $$\left(R. \alpha^{top}\right)_{g,n} = C_{g,n}(r,s).$$
\end{proposition}

\section{From the spectral curve to the Givental R-matrix}\label{sec:CurvetoGivental}

In this section we recall the main result of~\cite{DOSS,E}, which expresses the correlation differentials $\omega_{g,n}$ of the CEO topological recursion in terms of integral over moduli spaces of curves (Theorem \ref{thm:CohFT}). We then recall the identification \cite{DOSS} between topological recursion and Givental theory and apply it to a particular spectral curve (Definition \ref{def:Srs}). The result is an explicit expression for the coefficients of $\omega_{g,n}$ (Theorem \ref{thm:srequivalence}). 
\subsection{Local topological recursion}
The local version of the CEO topological recursion takes as input the following set of data $\S = (\Sigma, x, y, B)$:
\begin{enumerate}
\item[I).] A local spectral curve $\Sigma=\sqcup_{i=1}^r U_i$, given by the disjoint union of open disks with the center points $p_i$, $i=1,\dots,r.$
\item[II).] A holomorphic function $x\colon \Sigma\to\bbC$ such that the zeros of its differential $dx$ are $p_1,\dots,p_r$. We will assume the zeroes of $dx$ to be simple.
\item[III).] A holomorphic function $y\colon \Sigma\to\bbC$ which does not vanish on the zeroes of $dx$.
\item[IV).] A symmetric bidifferential $B$ defined on $\Sigma\times \Sigma$ with a double pole on the diagonal with residue $1$. 
\end{enumerate}

The output of the topological recursion procedure consists of a collection of symmetric differentials $\omega^{\S}_{g,n}$ defined on the topological product of the curve $\Sigma^{\times n}$. These correlation differentials take the following form:

\begin{theorem}\label{thm:CohFT} \cite{E,DOSS} The correlation differentials $\omega^{\S}_{g,n}$ produced via the topological recursion procedure from the specral curve $\S = (\Sigma,x,y,B)$ are equal to
\begin{equation*}\label{eq:C-fullformula}
 C^{2g-2+n}\!\!\!\! \sum_{\substack{i_1,\dots,i_n \\ d_1,\dots,d_n}}
\!\! \int_{\oM_{g,n}} \!\!\!\! \left(\Class^{\S}\right)_{g,n}(e_{i_1}\otimes \cdots \otimes e_{i_n}) 
\prod_{j=1}^n \psi_j^{d_j} \dd \left(\left( - \frac{1}{w_j}\frac{\dd}{\dd w_j}\right)^{d_j} \xi^{\S}_{i_j} \right).
\end{equation*}
\end{theorem}

\begin{remark}
The $\Class^{\S}$ defines a semi-simple CohFT, possibly with a non-flat unit. In this paper we will restrict the attention to CohFT with flat unit and we will write
$$\left(\Class^{\S}\right)_{g,n} = \left(R^{\S}.\alpha^{\S, top}\right)_{g,n}$$
to indicate its Givental decomposition (see Section \ref{sec:giv}; for CohFT with non-flat unit see \cite{PPZ} or \cite{LPSZ}, Section 2.3). 
\end{remark}
Let us describe the ingredients in the formula above in terms of the data of the specral curve, following \cite{DOSS}. The only difference with the usual representation is that we incorporate a torus action on cohomological field theories, fixing a point $(C,C_1,\dots,C_r)\in (\bbC^*)^{r+1}$. This formula doesn't depend on these parameters, though all its ingredients do.
\begin{itemize}
\item[i).] The local coordinates $w_i$ on $U_i$, $i=1,\dots,r$, are chosen such that $w_i(p_i)=0 $ and $ x=(C_iw_i)^2+x_i$

\item[ii).]The underlying topological field theory is given in idempotent basis by 
\begin{align*}\label{eq:C-underlyingTFT}
& \eta(e_i,e_j)= \delta_{ij}, \\ \notag
& \alpha^{\S, top}_{g,n}(e_{i_1}\otimes \cdots \otimes e_{i_n}) = \delta_{i_1\dots i_n} \left(-2C_i^2 C \frac{\dd y}{\dd w_i}(0)\right)^{-2g+2-n}.
\end{align*}

\item[iii).]The Givental matrix $R^{\S}(\zeta)$ is given by
\begin{equation*}\label{eq:C-matrixR}
-\frac 1\zeta (R^{\S})^{-1}(\zeta)_i^j=\frac{1}{\sqrt{2\pi\zeta}} \int_{-\infty}^\infty \left. \frac{B(w_i,w_j)}{\dd w_i}\right|_{w_i=0} \!\!\!\!\!\!\!\! \cdot  e^{-\frac{w_j^2}{2\zeta}}.
\end{equation*}

\item[iv).]The auxiliary functions $\xi^{\S}_i \colon \Sigma\to\bbC$ are given by
\begin{equation*}\label{eq:C-functionsxi}
\xi^{\S}_i(x):=\int^x \left.\frac{B(w_i,w)}{\dd w_i}\right|_{w_i=0}
\end{equation*}

\item[v).] \textit{DOSS Test} (see \cite{DNOPS}, Section 4): The following condition for the function $y$ is necessary and sufficient in order for the unit of the cohomological field theory $R^{\S}.\alpha^{\S}$ to be flat.
\begin{equation*}\label{eq:condition-y}
\frac{2C_i^2 C}{\sqrt{2\pi\zeta}} \int_{-\infty}^\infty \dd y\cdot e^{-\frac{w_i^2}{2\zeta}} = \sum_{k=1}^r (R^{-1})^i_k \left(2C_k^2 C \frac{\dd y}{\dd w_k}(0)\right)
\end{equation*}
\end{itemize}

\subsection{The spectral curve $\S_{r,s}$ and its Givental R-matrix}
Let us apply the formula in Theorem \ref{thm:CohFT} to the specral curve introduced in \cite{LPSZ}.
\begin{definition}\label{def:Srs} For $0 \leq s \leq r$, with $r$ and $s$ integer parameters, let the spectral curve $\S_{r,s}$ be defined by the following initial data in terms of a global coordinate $z$ on the Riemann sphere:
\begin{equation*}\label{eq:spectral-curve-data}
	\S_{r,s} := \left( \CP1, x(z)=-z^r+\log z,
 y(z)=z^s,  B(z,z') = \frac{\dd z\, \dd z'}{(z-z')^2} \right)
\end{equation*}
\end{definition}


\begin{proposition}[\cite{LPSZ}] \label{lem:flat-basis}
 The ingredients described above for the spectral curve $\S_{r,s}$ are given by:

\begin{itemize}
\item[i).] Choose the constants $C_i = 1/\sqrt{-2r}$ for $i=1,\dots,r$, and $C = r^{1 + s/r}/s$. With this choice the local coordinates $w_i$ on $U_i$, $i=1,\dots,r$ satisfy $ x= -\frac{w_i^2}{2r}+x(p_i)$.

\item[ii).]The underlying topological field theory is given by 
\begin{align*}
 &\eta(v_a,v_b)=\frac{1}{r}\delta_{a+b\mod r}\\
 &\alpha^{\S_{r,s}, top}_{g,n}(v_{a_1}\otimes \cdots \otimes v_{a_n}) = r^{2g-1} \delta_{a_1+\cdots+a_n-s(2g-2+n) \mod r}.
\end{align*}
where the flat coordinates $v_a$ are defined in terms of the idempotents by $$v_a := \sum_{i=0}^{r-1} \frac{J^{ai}}{r} e_i, \qquad a = 1, \dots, r.$$

\item[iii).]The Givental matrix $R^{\S_{r,s}}(\zeta)$ is given by
\begin{equation*}
R^{\S_{r,s}}(\zeta)= \exp\left(-\sum_{k=1}^{\infty} \frac{\mathrm{diag}_{a=1}^{r} B_{k+1}\left(\frac a r\right)}{k(k+1)}\zeta^k\right).
\end{equation*}

\item[iv).]The auxiliary functions $\xi^{\S_{r,s}}_a \colon \Sigma\to\bbC$ are given by
\begin{equation*}
\xi^{\S_{r,s}}_a=r^{\frac{r-a}r} \sum_{n=0}^\infty \frac{(nr+r-a)^n}{n!}e^{(nr+r-a)x}
\end{equation*}
\item[v).] DOSS Test is satisfied.
\end{itemize}
\end{proposition}

Let us now put together the ingredients for the correlation differentials $\omega_{g,n}^{\S_{r,s}}$, as in Theorem \ref{thm:CohFT}. First of all, in the coordinates $v_i$, the TFT and the Givental $R$-matrix for $\S_{r,s}$ coincide with the ones resulting in the Chiodo CohFT for the same parameters $r$ and $s$ by Proposition \ref{prop:RChiodo}.
Secondly, since $\frac{1}{r} \frac{d}{dx} = -\frac{1}{w_i}\frac{d}{dw_i}$, the computation of the derivatives of the auxiliary functions, divided by the powers of $r$, read
\begin{align*} \dd \left(\left( - \frac{\psi_j}{w_j}\frac{\dd}{\dd w_j}\right)^{d_j} \!\!\! \xi^{\S_{r,s}}_{i_j} \right) & = \dd \left[\left(\frac{\psi_j}{r} \frac{\dd}{\dd x_j}\right)^{d_j}
	\sum_{l_j=0}^\infty
	 \frac{(l_jr+r-a_j)^n}{l_j!}e^{(l_jr+r-a_j)x_j}\right]\\
	 &= 	 \dd_j
	 	 \sum_{l_j=0}^\infty 
	 \frac{(l_jr+r-a_j)^{l_j+d_j}}{l_j!} \left(\frac{\psi_j}{r}\right)^{d_j}e^{(l_jr+r-a_j)x_j}\\
	 	 &= 	 \dd_j 
	 	 \sum_{\mu_j=1}^\infty
	 \frac{\mu_j^{[\mu_j]}}{[\mu_j]!}\left(\frac{\mu_j}{r} \psi_j\right)^{d_j}e^{\mu_jx_j},\\
\end{align*}
where we write the euclidean division by $r$ as $\mu = [\mu]r + \langle \mu \rangle$, with $ \langle \mu \rangle < r$. The powers of $s$ only come from $C^{2g - 2 +n}$, and hence they are equal to $-(2g - 2 + n)$. The remaining powers of $r$ to compute amount to 
$$2g - 2 + n + \frac{(2g - 2 + n)s + \sum_{i} (r - a_i)}{r} =
2g - 2 + n + \frac{(2g - 2 + n)s + \sum_{i} \langle \mu_i \rangle}{r} ,$$ though it is handy to collect an extra $r^{\sum_i [\mu_i]}$ outside the product. This proves:
\begin{theorem}[\cite{LPSZ}]\label{thm:srequivalence}
The correlation differentials $\omega_{g,n}^{\S_{r,s}}$ of the spectral curve~\eqref{eq:spectral-curve-data} are equal to
	\begin{align*}
\dd_1\otimes \dots \otimes \dd_n \frac{r^{2g-2+n+b}}{s^{2g-2+n}} \prod_{j=1}^{n} \frac{\left(\frac{\mu_j}{r}\right)^{[\mu_j]}}{[\mu_j]!}  \!\!\!\!\!\! \sum_{\mu_1,\dots,\mu_n=1}^\infty	 
	\int_{\oM_{g,n}} \!\!\! \frac{\C_{g,n} \left (r,s; r \!-\! \langle\vec{\mu} \rangle \right)}{\prod_{j=1}^n (1-\frac{\mu_i}{r}\psi_i)}
	\; e^{\sum \mu_j x_j}
	\end{align*}
	where $b(r,s) = \left((2g-2+n)s+\sum_{j=1}^n \mu_j\right)/r$.
\end{theorem}

\begin{remark}
	Note that the case $s=1$ reproduces Theorem 1.7 in~\cite{SSZ}.
\end{remark}

Expanding the correlation differentials as
\begin{equation}
  \omega^{\S_{r,s}}_{g,n} =  \dd_1\otimes \cdots \otimes \dd_n \sum_{\mu_1,\dots,\mu_n=1}^\infty
	\frac{N_{g,\vec{\mu}}^{\S_{r,s}}}{b(r,s)!} \ e^{\sum_{j=1}^n \mu_j x_j} ,
	\end{equation}
we find:
	\begin{corollary}[\cite{LPSZ}]\label{cor:rsN}
\begin{align}
N^{\S_{r,s}}_{g,\vec{\mu}} = 
  b(r,s)! \frac{r^{b(r,s) + 2g - 2 +n}}{s^{2g - 2 +n}}
  \prod_{i=1}^{n}\frac{\left(\frac{\mu_i}{r}\right)^{[\mu_i]}}{[\mu_i]!}
  \int_{\overline{\mathcal{M}}_{g,n}}
  \frac{C_{g,n}(r,s;r - \langle \vec{\mu} \rangle)}{\prod_{j=1}^{n} (1 - \frac{\mu_j}{r} {\psi}_j)}.
\end{align}
\end{corollary}

\section{Equivalence statements: a new proof of the Johnson-Pandharipande-Tseng formula}\label{sec:JPT}

In this section we consider the case $s=r$ of Theorem \ref{thm:srequivalence}. In this case, the correlation differentials of this spectral curve are known to give the $r$-orbifold Hurwitz numbers $h^{\circ, [r]}_{g;\vec{\mu}}$ \cite{ BHLM, DLN, DLPS}, which enumerate connected Hurwitz coverings of the $2$-sphere of degree $|\vec{\mu}|$ and genus $g$, where the partition $\vec{\mu}$ determines the ramification profile over zero, the ramification over infinity is of cycle type $(r, r, \dots, r)$ and all other ramifications are simple.

Corollary \ref{cor:rsN} for $s=r$ specialises to
\begin{corollary}[\cite{LPSZ}]\label{TH2}
\begin{align}
N^{\S_{r,r}}_{g,\vec{\mu}} = b(r,r)!
  r^{b(r,r) }
  \prod_{i=1}^{n}\frac{\left(\frac{\mu_i}{r}\right)^{[\mu_i]}}{[\mu_i]!}
  \int_{\overline{\mathcal{M}}_{g,n}}
  \frac{C_{g,n}(r,r;r - \langle \vec{\mu} \rangle)}{\prod_{j=1}^{n} (1 - \frac{\mu_j}{r} {\psi}_j)}.
\end{align}
\end{corollary}
Plugging the $r-$orbifold Hurwitz numbers and the curve $\S_{r,r}$ into the TR-ELSV equivelance statement in Definition \ref{def:eqst}, we get:
\begin{corollary}\label{cor:eqorbifold}
The two propositions are equivalent:
\begin{align*}
\text{i)}& \qquad  h_{g, \vec{\mu}}^{\circ, [r]} = N^{\S_{r,r}}_{g, \vec{\mu}}\\
 \text{ii)}& \qquad   h_{g, \vec{\mu}}^{\circ, [r]} = b(r,r)!
  r^{b(r,r) }
  \prod_{i=1}^{n}\frac{\left(\frac{\mu_i}{r}\right)^{[\mu_i]}}{[\mu_i]!}
  \int_{\overline{\mathcal{M}}_{g,n}}
  \frac{C_{g,n}(r,r;r - \langle \vec{\mu} \rangle)}{\prod_{j=1}^{n} (1 - \frac{\mu_j}{r} {\psi}_j)}.
\end{align*}
\end{corollary}

On the other hand, $r$-orbifold Hurwitz numbers are also known to satisfy the John\-son-Pandharipande-Tseng (JPT) ELSV-type formula~\cite{JPT} (specialised here to the case
${G} = \mathbb{Z}/r\mathbb{Z}$, $U$ equal to the representation that sends $1$ to $e^{\frac{2 \pi i}{r}}$, and empty $\gamma$):

\begin{align} \label{eq:particular-jpt}
 h_{g;\vec{\mu}}^{\circ, [r]} = 
  b(r,r)! r^{b(r,r)}
  \prod_{i=1}^{n}\frac{\left(\frac{\mu_i}{r}\right)^{[\mu_i]}}{[\mu_i]!}
  \int_{\overline{\mathcal{M}}_{g,n}}
  \frac{p_*\sum_{i\geq 0} (-1)^i \lambda_i}{\prod_{j=1}^{n} (1 - \frac{\mu_j}{r} {\psi}_j)},
\end{align}
\begin{remark}
Note that the powers of $r$ are here slightly rearranged to easily match the equation for the correlation differentials.
\end{remark}
The class $p_*\sum_{i\geq 0} (-r)^i \lambda_i$ is described in~\cite{JPT} via admissible covers, while Chiodo's classes rely on the moduli space of $r$-th tensor roots. These two approaches are in fact equivalent:
\begin{proposition}[\cite{LPSZ}]\label{prop:Dima}
$
p_*\sum_{i\geq 0} (-r)^i \lambda_i = \C_{g,n}(r,r;r- \langle \vec{\mu} \rangle).
$
\end{proposition}
By Proposition \ref{prop:Dima} and Formula \eqref{eq:particular-jpt}, we can re-state Corollary \ref{cor:eqorbifold} as: 

\begin{corollary}[\cite{LPSZ}]\label{cor:eqorbifold2}
The two statements are equivalent:
\begin{enumerate}
\item[i)] The $r$-orbifold Hurwitz numbers satisfy the topological recursion from the spectral curve $\S_{r,r}$.\\
 \item[ii)]  The JPT formula holds.
\end{enumerate}
\end{corollary}

Since the JPT formula is proved independently from the topological recursion, Corollary \ref{cor:eqorbifold2} provides a new proof of the topological recursion for the numbers $ h_{g;\vec{\mu}}^{\circ, [r]}$. \\
On the other hand, both \cite{ BHLM, DLN} derive the topological recursion for $h_{g;\vec{\mu}}^{\circ, [r]}$ combining the cut-and-join equation with a polynomiality property which is extracted from the JPT formula itself. Hence one cannot conclude that Corollary \ref{cor:eqorbifold2} provides a new proof of the the JPT formula, unless this polynomiality property is derived independently from JPT.
This polynomiality property is proved in \cite{DLPS} with no use of JPT, see also \cite{KLS}. Therefore these results together, and via Corollary \ref{cor:eqorbifold2}, provide a new proof of Johnson-Pandharipande-Tseng formula.
\vspace{0.5cm}


\begin{changemargin}{-1.5cm}{2cm} 
{\renewcommand{\arraystretch}{1.2} 

\label{my-label}
\begin{tabular}{| l | r c  p{8.2cm} | }
\bottomrule
Case 
&
 Topological recursion 
 & & 
 ELSV-type formula 
 \\
\toprule
\begin{tabular}{@{}c@{}}
$s=1$\\
$r=1$ \\
\end{tabular}
& 
\begin{tabular}{r@{}}
\textit{The standard Hurwitz }\\
\textit{numbers $h^{\circ}_{g, \vec{\mu}}$}\\
\textit{are generated by}\\
$
\begin{cases}[r]
-z+\log z = x(z)\\
z = y(z)
\end{cases}
$
\\
 \end{tabular}
& $\iff$ & 
\begin{tabular}{l@{}}
 \textit{The ELSV formula}\\
 \begin{minipage}{2cm}
 \tiny
\begin{align*}
 \frac{ h^{\circ}_{g;\vec{\mu}} }{b!} = 
  \prod_{i=1}^{n}\frac{\mu_i^{\mu_i}}{\mu_i !}
  \int_{\oM_{g,n}}
  \frac{\sum_{i=0}^g (-1)^i \lambda_i}{\prod_{j=1}^{n} (1 - \mu_j {\psi}_j)},\\
 \end{align*}
\end{minipage}\\
\textit{holds, where: }$b = 2g - 2 + n + |\vec{\mu}|$. \\
\end{tabular}
\\
\midrule
 $s=r$ 
 & 
\begin{tabular}{r@{}}
\textit{The $r$-orbifold Hurwitz}\\
\textit{ numbers $h^{\circ, r}_{g, \vec{\mu}}$}\\
\textit{are generated by}\\
$
\begin{cases}[r]
-z^r+\log z = x(z)\\
z^r = y(z)
\end{cases}
$
\\
 \end{tabular}
 &  $\iff$ &  
\begin{tabular}{l@{}}
 \textit{The Johnson-Pandharipande-Tseng formula}\\
 \begin{minipage}{2cm}
 \tiny
\begin{align*}
 \frac{ h^{\circ,[r]}_{g;\vec{\mu}} } {b!} = 
  r^b
 \prod_i\frac{\left(\frac{\mu_i}{r}\right)^{[\mu_i]}}{[\mu_i]!}
  \int_{\overline{\mathcal{M}}_{g,n}}
  \frac{p_*\sum_{i\geq 0}^g (-1)^i \lambda_i}{\prod_{j=1}^{n} (1 - \frac{\mu_j}{r} {\psi}_j)},\\
 \end{align*}
\end{minipage}\\
\textit{holds, where: }$b = 2g - 2 + n + |\vec{\mu}|/r$. \\
\end{tabular}
 \\
 \midrule
 $s=1$ 
 &
\begin{tabular}{r@{}}
\textit{The $r$-spin Hurwitz}\\
\textit{ numbers $h^{\circ, r-spin}_{g, \vec{\mu}}$}\\
\textit{are generated by}\\
$
\begin{cases}[r]
-z^r+\log z = x(z)\\
z = y(z)
\end{cases}
$
\\
 \end{tabular}
 &  $\iff$ & 
\begin{tabular}{l@{}}
 \textit{The $r$-spin ELSV formula}\\
 \begin{minipage}{2cm}
 \tiny
\begin{align*}
 \frac{ h_{g;\vec{\mu}}^{\circ, r-spin}} {b!} = 
  r^{b-\chi}
  \prod_i\frac{\left(\frac{\mu_i}{r}\right)^{[\mu_i]}}{[\mu_i]!}
  \int_{\overline{\mathcal{M}}_{g,n}}
 \!\!\! \frac{\C_{g,n}(r,1, r - \langle\vec{\mu}\rangle)}{\prod_{j=1}^{n} (1 - \frac{\mu_j}{r} {\psi}_j)},\\
 \end{align*}
\end{minipage}\\
\textit{holds, where: }$b = (2g - 2 + n + |\vec{\mu}|)/r$. \\
\end{tabular}
\\
 \midrule
\end{tabular}
}
\end{changemargin}

\begin{remark}
In the table all the spectral curves are of genus zero, i.e. $\Sigma =  \mathbb{C} \mathbb{P}^1$, and hence with the standard kernel $B(z,z') = \frac{\dd z\, \dd z'}{(z-z')^2}$.
\end{remark}

\end{document}